%%%%%%%%%%%%%%%%%%%%%%%%%%%%%%% Macro Definitions %%%%%%%%%%%%%%%%%%%%%%%%%%%

\def\title#1{{\titlefont\noindent #1\bigskip}}

\def\author#1{{\largefont\noindent #1}\medskip}

\def\beginlinemode{\endmode
 \begingroup\obeylines\def\endmode{\par\endgroup}}
\let\endmode=\par

\newbox\theaddress
\def\address{\smallskip\beginlinemode\parindent 0in\getaddress}
{\obeylines
\gdef\getaddress #1 
 #2
 {#1\gdef\addressee{#2}%
   \global\setbox\theaddress=\vbox\bgroup\raggedright%
    \everypar{\hangindent2em}#2
   \def\endaddress{\egroup\endgroup \copy\theaddress \medskip}}}

\def\thanks#1{\footnote{}{\eightpoint #1}}

\long\def\Abstract#1{{\eightpoint\narrower\vskip\baselineskip\noindent
#1\smallskip}}

\def\skipfirstword#1 {}

\def\ir#1{\csname #1\endcsname}

\newdimen\currentht
\newbox\droppedletter
\newdimen\droppedletterwdth
\newdimen\drophtinpts
\newdimen\dropindent

\def\irrnSection#1#2{\edef\tttempcs{\ir{#2}}
\vskip6pt\penalty-3000
{\bf\noindent \expandafter\skipfirstword\tttempcs. #1}
\vskip6pt}

\def\irSubsection#1#2{\edef\tttempcs{\ir{#2}}
\vskip\baselineskip\penalty-3000
{\bf\noindent \expandafter\skipfirstword\tttempcs. #1}
\vskip6pt}

\def\irSubsubsection#1#2{\edef\tttempcs{\ir{#2}}
\vskip\baselineskip\penalty-3000
{\bf\noindent \expandafter\skipfirstword\tttempcs. #1}
\vskip6pt}

\def\References{\vbox to.25in{\vfil}\noindent{}{\bf References}
\vskip6pt\par}

\def\References{\vskip6pt\noindent{}{\bf References}
\vskip6pt\par}

\def\baselinebreak{\par \ifdim\lastskip<6pt
         \removelastskip\penalty-200\vskip6pt\fi}

\long\def\prclm#1#2#3{\baselinebreak
\noindent{\bf \csname #2\endcsname}:\enspace{\sl #3\par}\baselinebreak}

\def\Prf{\noindent{\bf Proof}: }

\def\rem#1#2{\baselinebreak\noindent{\bf \csname #2\endcsname}:\enspace }

\def\qed{{$\diamondsuit$}\vskip6pt}

\def\bibitem#1{\par\indent\llap{\rlap{\bf [#1]}\indent}\indent\hangindent
2\parindent\ignorespaces}

\long\def\eatit#1{}

\def\leftheadlinetext{}
\def\rightheadlinetext{}

\def\leftheadline{{\eightrm\folio\hfil \leftheadlinetext\hfil}}
\def\rightheadline{{\eightrm\hfil\rightheadlinetext\hfil\folio}}

\headline={\ifnum\pageno=1\hfil\else
\ifodd\pageno\rightheadline\else\leftheadline\fi\fi}

\def\tenpoint{\def\rm{\fam0\tenrm}
\textfont0=\tenrm \scriptfont0=\sevenrm \scriptscriptfont0=\fiverm
\textfont1=\teni \scriptfont1=\seveni \scriptscriptfont1=\fivei
\def\mit{\fam1} \def\oldstyle{\fam1\teni}
\textfont2=\tensy \scriptfont2=\sevensy \scriptscriptfont2=\fivesy
\def\cal{\fam2}
\textfont3=\tenex \scriptfont3=\tenex \scriptscriptfont3=\tenex
\def\it{\fam\itfam\tenit} % \it is family 4
\textfont\itfam=\tenit
\def\sl{\fam\slfam\tensl} % \sl is family 5
\textfont\slfam=\tensl
\def\bf{\fam\bffam\tenbf} % \bf is family 6
\textfont\bffam=\tenbf \scriptfont\bffam=\sevenbf
\scriptscriptfont\bffam=\fivebf
\def\tt{\fam\ttfam\tentt} % \tt is family 7
\textfont\ttfam=\tentt
\normalbaselineskip=12pt
\setbox\strutbox=\hbox{\vrule height8.5pt depth3.5pt  width0pt}%
\normalbaselines\rm}

\def\eightpoint{\def\rm{\fam0\eightrm}%
\textfont0=\eightrm \scriptfont0=\sixrm \scriptscriptfont0=\fiverm
\textfont1=\eighti \scriptfont1=\sixi \scriptscriptfont1=\fivei
\def\mit{\fam1} \def\oldstyle{\fam1\eighti}%
\textfont2=\eightsy \scriptfont2=\sixsy \scriptscriptfont2=\fivesy
\def\cal{\fam2}%
\textfont3=\tenex \scriptfont3=\tenex \scriptscriptfont3=\tenex
\def\it{\fam\itfam\eightit} % \it is family 4
\textfont\itfam=\eightit
\def\sl{\fam\slfam\eightsl} % \sl is family 5
\textfont\slfam=\eightsl
\def\bf{\fam\bffam\eightbf} % \bf is family 6
\textfont\bffam=\eightbf \scriptfont\bffam=\sixbf
\scriptscriptfont\bffam=\fivebf
\def\tt{\fam\ttfam\eighttt} % \tt is family 7
\textfont\ttfam=\eighttt
\normalbaselineskip=9pt%
\setbox\strutbox=\hbox{\vrule height7pt depth2pt  width0pt}%
\normalbaselines\rm}

\def\largefont{\def\rm{\fam0\largerm}
\textfont0=\largerm \scriptfont0=\largescriptrm \scriptscriptfont0=\tenrm
\textfont1=\largei \scriptfont1=\largescripti \scriptscriptfont1=\teni
\def\mit{\fam1} \def\oldstyle{\fam1\teni}
\textfont2=\largesy %\scriptfont2=\sevensy \scriptscriptfont2=\fivesy
\def\cal{\fam2}
%\textfont3=\largeex %\scriptfont3=\tenex \scriptscriptfont3=\tenex
\def\it{\fam\itfam\largeit} % \it is family 4
\textfont\itfam=\largeit
\def\sl{\fam\slfam\largesl} % \sl is family 5
\textfont\slfam=\largesl
\def\bf{\fam\bffam\largebf} % \bf is family 6
\textfont\bffam=\largebf %\scriptfont\bffam=\sevenbf \scriptscriptfont\bffam=\fivebf
\def\tt{\fam\ttfam\largett} % \tt is family 7
\textfont\ttfam=\largett
\normalbaselineskip=17.28pt
\setbox\strutbox=\hbox{\vrule height12.25pt depth5pt  width0pt}%
\normalbaselines\rm}

\def\titlefont{\def\rm{\fam0\titlerm}
\textfont0=\titlerm \scriptfont0=\largescriptrm \scriptscriptfont0=\tenrm
\textfont1=\titlei \scriptfont1=\largescripti \scriptscriptfont1=\teni
\def\mit{\fam1} \def\oldstyle{\fam1\teni}
\textfont2=\titlesy %\scriptfont2=\sevensy \scriptscriptfont2=\fivesy
\def\cal{\fam2}
%\textfont3=\largeex %\scriptfont3=\tenex \scriptscriptfont3=\tenex
\def\it{\fam\itfam\titleit} % \it is family 4
\textfont\itfam=\titleit
\def\sl{\fam\slfam\titlesl} % \sl is family 5
\textfont\slfam=\titlesl
\def\bf{\fam\bffam\titlebf} % \bf is family 6
\textfont\bffam=\titlebf %\scriptfont\bffam=\sevenbf \scriptscriptfont\bffam=\fivebf
\def\tt{\fam\ttfam\titlett} % \tt is family 7
\textfont\ttfam=\titlett
\normalbaselineskip=24.8832pt
\setbox\strutbox=\hbox{\vrule height12.25pt depth5pt  width0pt}%
\normalbaselines\rm}

\nopagenumbers

\font\eightrm=cmr8
\font\eighti=cmmi8
\font\eightsy=cmsy8
\font\eightbf=cmbx8
\font\eighttt=cmtt8
\font\eightit=cmti8
\font\eightsl=cmsl8
\font\sixrm=cmr6
\font\sixi=cmmi6
\font\sixsy=cmsy6
\font\sixbf=cmbx6

\font\largerm=cmr12 at 17.28pt
\font\largei=cmmi12 at 17.28pt
\font\largescriptrm=cmr12 at 14.4pt
\font\largescripti=cmmi12 at 14.4pt
\font\largesy=cmsy10 at 17.28pt
\font\largebf=cmbx12 at 17.28pt
\font\largett=cmtt12 at 17.28pt
\font\largeit=cmti12 at 17.28pt
\font\largesl=cmsl12 at 17.28pt

\font\titlerm=cmr12 at 24.8832pt
\font\titlei=cmmi12 at 24.8832pt
\font\titlesy=cmsy10 at 24.8832pt
\font\titlebf=cmbx12 at 24.8832pt
\font\titlett=cmtt12 at 24.8832pt
\font\titleit=cmti12 at 24.8832pt
\font\titlesl=cmsl12 at 24.8832pt

\tenpoint

%%%%%%%%%%%%%%%%%%%%%%%%%%%%%%% Internal References %%%%%%%%%%%%%%%%%%%%%%%%%%%

%%%%%%%%%%%%%%%%%%%%%%%%%%%%%%% Begin Paper %%%%%%%%%%%%%%%%%%%%%%%%%%%

%\input macros19Jul95

\def\manyby{\hbox to.75in{\hrulefill}}
\hsize 6.5in 
\vsize 9.2in

\tolerance 3000
\hbadness 3000

\def\refCam{Cam}
\def\refCat{Cat}
\def\refCMa{CM1}
\def\refCMb{CM2}
\def\refDGM{DGM}
\def\refDu{Du}
\def\refFione{F1}
\def\refFitwo{F2}
\def\refFHH{FHH}
\def\refGGR{GGR}
\def\refGi{Gi}
\def\refvanc{Ha1}
\def\refsundance{Ha2}
\def\refravello{Ha3}
\def\refanti{Ha4}
\def\reffatpts{Ha5}
\def\refigp{Ha6}
\def\refHi{Hi}
\def\refId{Id}
\def\refMir{Mi}
\def\refMu{Mu}
\def\refNone{N1}
\def\refNtwo{N2}
\def\refR{R}

\def\item#1{\par\indent\indent\llap{\rlap{#1}\indent}\hangindent
2\parindent\ignorespaces}

\def\itemitem#1{\par\indent\indent
\indent\llap{\rlap{#1}\indent}\hangindent
3\parindent\ignorespaces}

\def\C#1{\hbox{$\cal #1$}}

\def\div#1 #2 #3 #4 #5 #6 #7 #8 {#1L #2E_1 #3E_2 #4E_3 #5E_4 #6E_5 #7E_6 #8E_7}

\def\pr#1{\hbox{{\bf P}${}^{#1}$}}
\def\leftheadlinetext{Harbourne, Holay and Fitchett}
\def\rightheadlinetext{Resolutions}

{\largefont
\hbox to\hsize{Resolutions of Ideals of Uniform Fat Point Subschemes of \pr2\hfil}
}
\vskip\baselineskip

{\font\authorrm=cmr12 at 14.4pt
\authorrm
\hbox to\hsize{Brian Harbourne\hfil }
}
\vskip6pt

\address
Department of Mathematics and Statistics, University of Nebraska-Lincoln
Lincoln, NE 68588-0323
email: bharbour@math.unl.edu
WEB: http://www.math.unl.edu/$\sim$bharbour/
\endaddress

{\font\authorrm=cmr12 at 14.4pt
\authorrm
\hbox to\hsize{Sandeep Holay\hfil }
}
\vskip6pt

\address
Department of Mathematics, Southeast Community College, Lincoln, NE 
email: sholay@math.unl.edu
\endaddress

{\font\authorrm=cmr12 at 14.4pt
\authorrm
\hbox to\hsize{Stephanie Fitchett\hfil }
}
\vskip6pt

\address
Department of Mathematics, Duke University, Durham, NC
email: sfitchet@math.duke.edu
\smallskip
June 4, 1999\endaddress
\vskip-\baselineskip

\thanks{\vskip -6pt
\noindent This work benefitted from a National Science Foundation grant.
\smallskip
\noindent 1991 {\it Mathematics Subject Classification. } 
Primary 13P10, 14C99. 
Secondary 13D02, 13H15.
\smallskip
\noindent {\it Key words and phrases. }  Ideal generation 
conjecture, symbolic powers,
resolution, fat points, maximal rank.\smallskip}

\vskip\baselineskip
\Abstract{Abstract: Let $I$ be the ideal
defining a set of general points $p_1,\ldots,p_n\in\pr2$.
There recently has been progress in showing that
a naive lower bound for the Hilbert functions
of symbolic powers $I^{(m)}$ is in fact attained
when $n>9$. Here, for $m$ sufficiently large, the minimal 
free graded resolution of $I^{(m)}$ is determined when $n>9$ 
is an even square, assuming only that this lower bound 
on the Hilbert function is attained.
Under ostensibly stronger conditions (that are nonetheless
expected always to hold), a similar result is shown to
hold for odd squares, and for infinitely many $m$ for 
each nonsquare $n$ bigger than 9. All results hold for
an arbitrary algebraically closed field $k$.}
\vskip\baselineskip

\irrnSection{Introduction}{intro}
The concerns of this paper, in which we 
fix an algebraically closed ground field $k$, are rooted in work by
Dubreil on numbers of generators of homogeneous ideals
and in work by Nagata resolving Hilbert's 14th Problem
and posing a still open conjecture regarding the minimum degree
of a curve with certain assigned multiplicities.

Let $I$ be a homogeneous ideal defined by base point conditions in \pr2
(more precisely, in terms defined below, let $I$ be the 
saturated homogeneous ideal
defining the fat point subscheme $m_1p_1+\cdots+m_np_n$
for points $p_i\in\pr2$ and nonnegative multiplicities $m_i$).
If the points are general and if $n>9$, Nagata [\refNone, \refNtwo] 
conjectures that the least degree of a nontrivial element of $I$ is 
always more than $(m_1+\cdots+m_n)/\sqrt{n}$. He proves this 
when $n$ is a square, which he then uses in his 
famous counterexample to Hilbert's 14th Problem. 

Nagata's work can be seen as giving
a bound on the Hilbert function of $I$ for low degrees.
Dubreil's work [\refDu], with additional developments
by [\refDGM] and [\refCam], 
derives from the Hilbert function of $I$ bounds on numbers of 
generators of $I$. Recently, by paying attention to geometry
(in particular, by keeping in mind the multiplicities $m_i$ attached 
to the points), [\refFione] improves these bounds under certain conditions.
In various special cases, the minimal free resolution for $I$
has now also been determined:
[\refCat] and [\ir{reffatpts}] do so for points on plane conics;
[\refFione, \refFitwo] does so for $n\le6$ general points, 
extended in unpublished work of the authors to $n=8$; 
[\refigp] does so for uniform subschemes $mp_1+\cdots+mp_n$ 
with $p_i$ general and $n\le 9$;
and [\refId] does so for $n>9$ for subschemes $2p_1+\cdots+2p_n$ when the points 
$p_i$ are general. 

In each of the cases that a resolution was determined, the Hilbert function of $I$
was known previously. One stumbling block to finding minimal free resolutions generally
is that the Hilbert function of I is not always known,
as evidenced by the fact that Nagata's conjecture is still open. However,
for subschemes $m_1p_1+\cdots+m_np_n$ with $n>9$ and $p_i$ general, evidence is
accumulating [\ir{refCMa}, \ir{refCMb}, \ir{refR}] that the Hilbert function 
takes a certain minimal form
if the multiplicities $m_1,\ldots,m_n$ are fairly uniform.
The conjecture that this is indeed what happens is formulated here precisely 
as \ir{qhc}. (In fact, Nagata's conjecture is a consequence  
of \ir{qhc}.)

Under the assumption that this minimal Hilbert function 
is under certain circumstances achieved, the main result of this paper, \ir{main}, 
shows for each $n>9$ that there are infinitely many $m$ 
such that the minimal free resolution for $I$ 
also takes a certain minimal form, when $I$ corresponds
to the subscheme $mp_1+\cdots+mp_n$ with $p_i$ general.
This suggests that this minimal form of the minimal free resolution 
may be attained for all $m$ for every $n>9$ when the points $p_i$ are general.
This suggestion is formalized by \ir{urc}. Additional evidence for
\ir{urc} is given in \ir{egs}.

Although our results are for \pr2, the general problem is of interest
for all \pr N, $N>1$. So let $p_1,\ldots,p_n\in\pr N$ 
be general points of projective space.
Let $P_j$ be the homogeneous ideal defining $p_j$
in the homogeneous coordinate ring $R=k[x_0,\ldots,x_N]$
of \pr N. An ideal of the form $I({\bf m};n,N)=P_1^{m_1}\cap \cdots\cap P_n^{m_n}$,
where ${\bf m}$ denotes the $n$-tuple 
$(m_1,\ldots,m_n)$ of nonnegative integers $m_j$,
defines a {\it fat point\/} subscheme $m_1p_1+\cdots+m_np_n$
of \pr N. (Thus, for example, $I((m,\ldots,m);n, N)$ is the $m$-th
symbolic power $I((1,\ldots,1);n,N)^{(m)}$ of the ideal $I((1,\ldots,1);n,N)$ generated by all
forms that vanish at the $n$ points.) Being a homogeneous ideal,
$I({\bf m};n,N)$ is a direct sum of its components $I({\bf m};n,N)_t$ of 
each degree $t$. We can regard $I({\bf m};n,N)_t$ as the linear system
of forms of degree $t$ with the imposed base point conditions
of vanishing at each point $p_j$ to order at
least $m_j$. 

A natural problem is that of determining the ($k$ vector space) dimension 
$h_{I({\bf m};n,N)}(t)$ of $I({\bf m};n,N)_t$
as a function of $t$; $h_{I({\bf m};n,N)}$ is known as the Hilbert function 
of $I({\bf m};n,N)$. A more difficult problem is to determine the minimal number 
of homogeneous generators of $I({\bf m};n,N)$ needed in each degree $t$;
this is equivalent to determining the Hilbert function of 
$I({\bf m};n,N)\otimes_R k$, since the dimension of the homogeneous
component $(I({\bf m};n,N)\otimes_R k)_t$ is just the 
minimal number of homogeneous
generators of $I({\bf m};n,N)$ needed in degree $t$. The gold standard in
connection to problems like these is the determination
of a minimal free resolution 
$0\to F_{N-1}\to \cdots\to F_0\to I({\bf m};n,N)\to 0,$
or at least a determination of the modules $F_j$,
up to isomorphism as graded $R$-modules,
since knowing the modules $F_j$ up to graded isomorphism
is sufficient to determine both the numbers of homogeneous
generators of $I({\bf m};n,N)$ in each degree and the Hilbert function
of $I({\bf m};n,N)$.

The situation for $N=1$ is trivial, since then saturated
homogeneous ideals are principal. For $N>1$, it is only for
$N=2$ that a general conjecture has been posed,
and then only for the Hilbert function
(see \ir{genconj} below, which uses the language of
linear systems of sections on blow ups of \pr2 at the
$n$ points). It is perhaps worth mentioning
that complete answers for $N=2$ have been attained for all of the 
questions raised above for $n\le 8$. Under different terminology,
the Hilbert function of the ideal $I$ corresponding to any
subscheme $m_1p_1+\cdots+m_np_n\in\pr2$ with $p_i$ general and $n\le 8$
was determined by Nagata [\ir{refNtwo}]. (More generally,
the Hilbert function of $I$ is known as long as a cubic curve passes
through the points $p_i$; see [\ir{refanti}].) 
For $n\le 8$, recent unpublished work [\ir{refFHH}]
of the authors has now determined the 
minimal free resolution of $I({\bf m};n,2)$, 
but the general situation remains unclear for $n>8$.

\irrnSection{Background}{bkgnd}
Since hereafter we restrict our attention to $N=2$, $R$ will denote
the homogeneous coordinate ring $k[x,y,z]$ of \pr2 and we will write
$I({\bf m};n)$ for $I({\bf m};n,2)$. Moreover, unless otherwise explicitly
mentioned, we will hereafter always assume 
that $I({\bf m};n)$ involves $n$ {\it general\/} 
points of \pr2.

So suppose one knows the Hilbert function $h_{I({\bf m};n)}$ of $I({\bf m};n)$.
The problem of then determining the number of generators
in each degree of a minimal set of homogeneous generators for $I({\bf m};n)$
is equivalent to the problem of determining the modules
in the resolution up to graded isomorphism.
This is because, since the subscheme $Z_{{\bf m},n}\subset\pr2$ 
defined by $I({\bf m};n)$ is arithmetically 
Cohen-Macaulay, the resolution takes the form
$0\to F_1\to F_0\to I({\bf m};n)\to 0$.
Knowing the number of generators in each degree
gives $F_0$, and, given the Hilbert function of $I({\bf m};n)$,
exactness of the sequence $0\to F_1\to F_0\to I({\bf m};n)\to 0$
allows one to then determine the Hilbert function of
$F_1$ and thus $F_1$ itself.

For $n\ge 9$, however, the solution to these equivalent problems,
of determining resolutions or of determining numbers of generators
for $I({\bf m};n)$ given its Hilbert function, is not yet clear even 
conjecturally. Things improve somewhat if we restrict to uniform ${\bf m}$:

\rem{Definition}{qunifdef} We will say that ${\bf m}$ or 
more generally $I({\bf m};n)$
is {\it uniform\/} if ${\bf m}=(m,\ldots,m)$ for some $m$ and
{\it quasiuniform\/} if ${\bf m}=(m_1,m_2,\ldots,m_r)$ with
$r\ge 9$ and $m_1=\cdots=m_{9}\ge m_{10}\ge \cdots \ge m_r$.
\vskip6pt

\rem{Notation}{primenot} It will be helpful to denote 
$I({\bf m},n)$ by $I(m;n)$ 
when the $n$-tuple ${\bf m}$ is $(m,\ldots,m)$, or more generally
by $I(m_1,m_2,\ldots, m_r;n)$ when ${\bf m}=(m_1,m_2,\ldots,m_r,\ldots,m_r)$.
Given ${\bf m}=(m_1,\ldots,m_n)$, we will denote
$(m_1+1,m_2\ldots,m_n)$ by ${\bf m}'$ and 
$(m_1-1,m_2\ldots,m_n)$ by ${\bf m}''$; likewise, we will write
$I({\bf m},n)'$ for $I({\bf m}',n)$ and 
$I({\bf m},n)''$ for $I({\bf m}'',n)$.  
\vskip6pt

All is known for $I(m;n)$ for $n\le9$ ([\ir{refigp}]), while for $n\ge 10$
we have \ir{qhc} (which is a special case of \ir{genconj}) for the Hilbert function,
and we have \ir{urc} for the resolution.

In preparation for stating these
conjectures in the next section, consider
the space $R_d$ of all forms of a given degree $d$. The subspace of
those vanishing at a given point with multiplicity $m$ or more
has dimension exactly $\hbox{max}\{0, ((d+1)(d+2)-m(m+1))/2\}$.
Thus a natural first conjecture is that the subspace of
$R_d$ of those vanishing with multiplicity $m_i$ or more
at each of $n$ general points $p_i$ has dimension 
$\hbox{max}\{0, ((d+1)(d+2)-\sum_{i=1}^nm_i(m_i+1))/2\}$. 
Accordingly, we make the following definition.

\rem{Definition}{unhindered} Given ${\bf m}=
(m_1,\ldots,m_n)$ with each $m_i\ge 0$, we will
say that $I({\bf m};n)$ is {\it unhindered\/} if $h_{I({\bf m};n)}(t)=
\hbox{max}\{0, ((t+1)(t+2)-\sum_{i=1}^nm_i(m_i+1))/2\}$ for all $t\ge 0$.
\vskip6pt

Sheafifying $I({\bf m};n)$ gives an ideal sheaf \C I,
and we note that $I({\bf m};n)$ being unhindered is the same thing
as $h^1(\pr2, \C I(t))$ vanishing for every $t$ for which 
$h^0(\pr2,\C I(t))>0$, where $\C I(t)$ denotes 
the twist $\C I\otimes \C O_{\pr2}(t)$ by $t$ times the class of a line.

Now, although $I({\bf m};n)$ is not in general unhindered, 
the known hindrances (which
\ir{genconj} is formulated to account for) are fairly special,
and in fact \ir{genconj} implies that all
quasiuniform $I({\bf m};n)$ are unhindered if $n\ge 10$
(see \ir{qunifunhin}).
Moreover, [\ir{refCMa}, \ir{refCMb}] shows that 
$I(m;n)$ is unhindered for all $n\ge 10$ with uniform
multiplicity $m\le 12$. Thus it is reasonable to
study resolutions of uniform ideals
$I(m;n)$ under assumptions of unhinderedness.

Note that if $I({\bf m};n)$ is unhindered,
it follows that the Castelnuovo-Mumford regularity of
$I({\bf m};n)$ is at most $\alpha({\bf m};n)+1$, where
$\alpha({\bf m};n)$ is the least degree $t$ such that $I({\bf m};n)_t\ne0$
(and where $\alpha(m;n)$ denotes the least degree $t$ such that $I(m;n)_t\ne 0$).
This implies (see [\ir{refDGM}] or Lemma 2.9 of [\reffatpts])
that the first syzygy module $F_0$ in a 
minimal free resolution of $I({\bf m};n)$ has generators
in at most two degrees, $\alpha({\bf m};n)$ and $\alpha({\bf m};n)+1$.
The number of generators in degree $\alpha({\bf m};n)$ is clearly
$h({\bf m};n)=h_{I({\bf m};n)}(\alpha({\bf m};n))$, while if
$I({\bf m};n)$ is unhindered the number of
generators in degree $\alpha({\bf m};n)+1$ is at least
$\hbox{max}\{0, \alpha({\bf m};n)+2-2h({\bf m};n)\}$, 
since $\alpha({\bf m};n)+2-2h({\bf m};n)$ is
the difference in dimensions of
$I({\bf m};n)_{\alpha({\bf m};n)+1}$ and 
$I({\bf m};n)_{\alpha({\bf m};n)}\otimes R_1$. Thus, if
$I({\bf m};n)$ is unhindered, the minimal
possible rank for $F_0$ is 
$\hbox{max}\{h({\bf m};n),\alpha({\bf m};n)+2-h({\bf m};n)\}$.
In this paper we will adduce circumstances in which
this minimum is achieved. Thus it is convenient to make
the following definition:

\rem{Definition}{rkmin} We will say that $I({\bf m};n)$
is {\it rank minimal\/} if $I({\bf m};n)$ is unhindered
and if the rank of $F_0$ in the minimal free
resolution $0\to F_1\to F_0\to I({\bf m};n)\to 0$ is 
$\hbox{max}\{h({\bf m};n),\alpha({\bf m};n)+2-h({\bf m};n)\}$.
\vskip6pt

An ideal $I({\bf m};n)$ being rank minimal  
is equivalent to $I({\bf m};n)$ being unhindered with the multiplication 
homomorphism $I({\bf m};n)_{\alpha({\bf m};n)}\otimes R_1\to 
I({\bf m};n)_{1+\alpha({\bf m};n)}$ having maximal rank (i.e,
either being injective or surjective). Moreover, 
an unhindered ideal $I({\bf m};n)$ being 
rank minimal is equivalent to the minimal free resolution taking
the following form (where $a=\alpha({\bf m};n)$, 
$h=h_{I({\bf m};n)}(a)$, $b=\hbox{max}\{0, a+2-2h\}$
and $c=\hbox{max}\{0, 2h-a-2\}$---hence either $b$ or $c$
is zero):

\vskip6pt
\noindent\hbox to\hsize{\hbox to0in{$(*)$\hss}\hfil
$0\to R[-a-2]^{a+1-h}\oplus R[-a-1]^{c}
\to R[-a-1]^{b}\oplus R[-a]^{h}\to I({\bf m};n)\to 0.$\hfil}
\vskip6pt

\noindent We recall that $R[a]^b$, for example, denotes 
the direct sum of $b$ copies of the homogeneous $R$-module $R[a]$,
with the grading defined by $R[a]_t=R_{t+a}$. 

To justify $(*)$, for an unhindered ideal note that a minimal free resolution
of the form $(*)$ immediately implies that $I({\bf m};n)$
is rank minimal. Conversely, if an unhindered ideal $I({\bf m};n)$
is rank minimal, then $F_0$ is as given
and one obtains $F_1$ as follows. 
Sheafifying (and suppressing ${\bf m}$ and $n$) gives
$0\to \C F_1\to \C F_0\to \C I\to 0$
and $0\to \C I\to \C O_{\pr 2}\to \C O_Z\to 0$,
where $\C F_0 = \C O_{\pr 2}(-a)^{\oplus h}
\oplus \C O_{\pr 2}(-a-1)^{\oplus b}$ and $Z$ is the 
subscheme defined by \C I.
Clearly, $h^2(\pr 2, \C F_0(a))$ vanishes,
while $h^1(\pr 2, \C I(a))=0$ since
$I({\bf m};n)$ is unhindered. From 
$0\to \C F_1(a)\to \C F_0(a)\to \C I(a)\to 0$ 
one sees that $h^0(\pr 2, \C F_1(a))=0$ and that $h^2(\pr 2, \C F_1(a))=0$
and hence that $\C F_1 = \C O_{\pr 2}(-a-1)^{\oplus c}
\oplus \C O_{\pr 2}(-a-2)^{\oplus d}$ for some $c$ and $d$.
To determine $c$ and $d$, take
$h^0$ of $0\to \C F_1(t)\to \C F_0(t)\to \C I(t)\to 0$
for $t=a+1$ and $t=a+2$.

Our main result can now be stated.

\prclm{Theorem}{main}{Consider $n>9$ general points $p_i\in \pr2$.
\item{(a)} Say $n$ is not a square. Then there are infinitely many
$m>0$ such that $I(m;n)$ is rank minimal
whenever $I(m;n)$ and $I(m;n)'$ are unhindered.
\item{(b)} Say $n$ is an odd square. Then for all $m\ge(n-9)/8$,
$I(m;n)$ is rank minimal whenever $I(m;n)$ and $I(m;n)'$ are unhindered.
\item{(c)} Say $n=r^2$ is an even square. Then for all $m\ge(r-2)/4$,
$I(m;n)$ is rank minimal whenever $I(m;n)$ is unhindered.}

\Prf Part (a) is \ir{nnsqr}, part (b) is handled by
\ir{evensqr}(b), and part (c) follows directly from 
\ir{evensqr}(a) and \ir{hlbfnctn}(i).
\qed

The method used here when $n$ is not an even square was employed in
[\ir{refigp}] for certain special values of $n$; thus
\ir{main} extends this to all $n$ which are not even squares. The
method of [\ir{refigp}] fails
when $n$ is an even square, so we develop a 
new approach special to that case.

The following result will be useful.
Define $q({\bf m};n)$ to be $h_{I({\bf m}';n)}(\alpha({\bf m};n))$,
and define $l({\bf m};n)$ to be 
$h_{I({\bf m}'';n)}(-1+\alpha({\bf m};n))$. We will also use $q(m;n)$
for $q({\bf m};n)$ and $l(m;n)$ for $l({\bf m};n)$ in case
${\bf m}=(m_1,\ldots,m_n)$ when $m=m_1=\cdots=m_n$.
%In addition, we define 
%$q^*({\bf m};n)=q({\bf m};n)-
%((\alpha({\bf m};n)+1)(\alpha({\bf m};n)+2)-(m_1+1)(m_i+2)-\sum_{i=2}^nm_i(m_i+1))/2$
%and $l^*({\bf m};n)=l({\bf m};n)-
%(\alpha({\bf m};n)(\alpha({\bf m};n)+1)-m_1(m_i+1)-\sum_{i=2}^nm_i(m_i+1))/2$,
%with $q^*(m;n)$ and $l^*(m;n)$ following the same convention as above.
%(We note that $q^*({\bf m};n)=h^1(\pr2,\C F(\alpha({\bf m};n)))$
%and $l^*({\bf m};n)=h^1(\pr2,\C G(\alpha({\bf m};n)-1))$, where \C F and \C G denote 
%the sheafifications of $I({\bf m}';n)$ and $I({\bf m}'';n)$, respectively.)

\prclm{Lemma}{igpLemma}{Assume $I({\bf m};n)$ is unhindered.
\item{(a)} If $q({\bf m};n)=0$ and
$l({\bf m};n)=0$, then $I({\bf m};n)$ is rank minimal.
\item{(b)} If ${\bf m}=(m_1,\cdots,m_n)$ has $m_1=m_2$, 
and $q({\bf m};n)=0$, then $I({\bf m};n)$ is rank minimal.
\item{(c)} If ${\bf m}=(m_1,\cdots,m_n)$ has $m_1=m_2$, 
$l({\bf m};n)>0$, and $I({\bf m};n)'$ and $I({\bf m};n)''$
are unhindered, then $I({\bf m};n)$ is rank minimal.
}

\Prf (a) Keeping in mind our comment about regularity 
(preceding \ir{rkmin}), this follows from 
Lemma 4.1 of [\ir{refigp}].

(b) Let $f$ be the linear form vanishing on the line through 
the points $p_1$ and $p_2$. Then multiplication by $f$ gives an injection
$I(-1+m_1,m_2,m_3,\cdots,m_n;n)_{-1+\alpha({\bf m};n)}
\to I(m_1,1+m_2,m_3,\cdots,m_n;n)_{\alpha({\bf m};n)}$, and since
the points $p_1,\cdots,p_n$ are general, we know
$I(m_1,1+m_2,m_3,\cdots,m_n;n)_{\alpha({\bf m};n)}$ and 
$I(1+m_1,m_2,\cdots,m_n;n)_{\alpha({\bf m};n)}$ have the same dimension,
hence $l({\bf m};n)\le q({\bf m};n)$. Therefore,
if $q({\bf m};n)=0$, then $I({\bf m};n)$ is rank minimal by (a).

(c) As in (b), we have $l({\bf m};n)\le q({\bf m};n)$,
but now $0<l({\bf m};n)$, so, applying the unhindered hypotheses,
the result follows from Lemma 4.2 of [\ir{refigp}]. \qed

\irrnSection{Conjectures}{cnjctrs}
Although \ir{main} employs certain assumptions about unhinderedness,
this is not a serious restriction since the assumptions are either known
to be met or, as we discuss below, are conjectured to be.  

To begin, the Hilbert function of $I({\bf m};n)$ 
is known (even if the points 
$p_i$ are not general) for $n\le9$ but it is somewhat complicated.
It remains unknown (but conjectured) for $n>9$, but, as indicated by
the following Quasiuniform Hilbert function Conjecture (QHC),
it is expected that the answer is fairly simple for quasiuniform ideals:

\prclm{Conjecture}{qhcWithAcronym}{For $n>9$ general points of
\pr2, $I({\bf m};n)$, $I({\bf m};n)'$ and $I({\bf m};n)''$
are unhindered if $I({\bf m};n)$ is quasiuniform.}

Similarly, whereas the problem of determining resolutions is 
known but complicated for ideals $I({\bf m};n)$ involving $n\le 8$ general points
(and for uniform ${\bf m}$ when $n=9$) and unknown and in general 
unconjectured for $n>9$, things
are expected to be particularly simple for uniform ideals
for $n>9$ general points. In particular, we have the following 
Uniform Resolution Conjecture (URC) (which is equivalent to 
Conjecture 6.3 (Maximal Rank Conjecture) of [\ir{refigp}]):

\prclm{Conjecture}{urcWithAcronym}{For $n>9$ general points of
\pr2, uniform unhindered ideals $I(m;n)$ are rank minimal.} 

\ir{egs} gives various examples in which URC is now known to hold.

Whereas it is unclear how URC might be extended to the nonuniform case,
QHC is, as we show in \ir{qunifunhin}, a special case of a more general
conjecture, \ir{genconj} given below, in which no assumptions are made 
on the uniformity of the multiplicities $m_i$. 
For other equivalent variants of
\ir{genconj}, see [\ir{refvanc}, \ir{refHi}, \ir{refGi}, \ir{refravello}];
for a nice survey, see [\refMir]. 

\prclm{Conjecture}{genconj}{For each positive integer $d$
there is a nonempty open set $U_d$ of points $(p_1,\ldots,p_n)\in(\pr2)^n$
such that on the surface $X$ obtained by blowing up 
$p_1,\ldots,p_n$ we have the following (where $L$ 
is the total transform of a line):
\item{(i)} either $h^0(X, \C O_X(F))=0$ or $h^1(X, O_X(F))=0$
for any numerically effective divisor $F$ 
with $F\cdot L\le d$; 
\item{(ii)} if $C$ is a prime divisor on $X$ of 
negative self-intersection with $C\cdot L\le d$,
then $C^2=C\cdot K_X=-1$.}

\rem{Remark}{algor} As an aside, we mention
how by assuming \ir{genconj} it is a simple matter to explicitly
compute the Hilbert function of an ideal $I({\bf m};n)$.
First,  $I({\bf m};n)_t$ corresponds for each $t$
to the complete linear system $|D_t|$ on $X$
of the divisor $D_t=tL-m_1E_1-\cdots-m_nE_n$,
where $X\to \pr2$ is the blow up of the points $p_1,\ldots,p_n$,
$L$ is the total transform to $X$ of a line in \pr2, and
each $E_i$ is the exceptional locus of the blow up of $p_i$.
Next, one subtracts off from $D_t$ irreducible effective divisors 
$C$ of the form $C^2=C\cdot K_X=-1$ (all such $C$ are known 
by [\ir{refNtwo}]) for which $C\cdot D_t<0$.
Eventually a divisor $D$ is obtained from $D_t$ which either meets all such $C$ 
nonnegatively or which has $D\cdot L<0$. In the latter case,
$|D|$ is empty and $h^0(X, \C O_X(D))=0$.
In the former case, \ir{genconj}(ii) implies either that
$|D|$ is empty or $D$ is numerically effective, but by duality $h^2(X, \C O_X(D))=0$ 
for any numerically effective divisor, so \ir{genconj}(i)
and Riemann-Roch allow one to compute $h^0(X, \C O_X(D))$
and hence the Hilbert function of $I({\bf m};n)$.
\vskip6pt

\rem{Remark}{qunifunhin} Using the method of proof of
Corollary 5.6 of [\ir{refigp}], we now show that \ir{genconj} implies
\ir{qhcWithAcronym}. Using the notation of \ir{algor},
let $D_t=tL-m_1E_1-\cdots-m_nE_n$, where by quasiuniformity
$m_1=\cdots=m_{9}$ (we will denote this common multiplicitiy by $m$),
and $m\ge m_{10}\ge \cdots$. But $C=3L-E_1-\cdots-E_{9}$ is
numerically effective (there is always an irreducible cubic
through 9 general points), so, if $t\ge \alpha((m_1,\ldots,m_n);n)$,
then $C\cdot D_t\ge 0$, hence $t\ge 3m$ so 
$D_t=-mK_X+(t-3m)L+(m-m_{10})E_{10}+\cdots+(m-m_n)E_n$. If $D_t$ were
not numerically effective, then there would be a reduced, irreducible
curve $A$ such that $A\cdot D_t<0$, but then we would have $A^2<0$
and $A\cdot L\le t$, hence, by \ir{genconj}, $A^2=A\cdot K_X=-1$.
Since $E_i\cdot D_t\ge 0$ for all $i$, we see 
$A$ is not $E_i$, but then 
$A\cdot D_t\ge A\cdot (-mK_X) = m\ge 0$, which shows that
$D_t$ is numerically effective if $h^0(X,\C O_X(D_t))>0$;
thus $I((m_1,\ldots,m_n);n)$ is unhindered by \ir{genconj}. To see that
$I((m_1,\ldots,m_n);n)'$ is unhindered, take $D_t$ to be
$tL-(m_1+1)E_1-m_2E_2-\cdots-m_nE_n$, where
$t\ge \alpha((m_1+1,m_2,\ldots,m_n);n)$; then $C\cdot D_t\ge 0$ implies 
$3t\ge 9m+1$ so $t\ge 3m+1$ so 
$D_t=-mK_X+(t-3m)L+(L-E_1)+(m-m_{10})E_{10}+\cdots+(m-m_n)E_n$
and the argument now proceeds as before.
Finally, for $I((m_1,\ldots,m_n);n)''$, take
$D_t=tL-(m_1-1)E_1-m_2E_2-\cdots-m_nE_n$, where
$t\ge \alpha((m_1-1,m_2,\ldots,m_n);n)$; then $3t\ge 9m-1$ so 
$t\ge 3m$ so 
$D_t=-mK_X+(t-3m)L+E_1+(m-m_{10})E_{11}+\cdots+(m-m_n)E_n$
and again we obtain the result.

\irrnSection{When $n$ is not a square}{ntvnsqrs}

\prclm{Proposition}{nnsqr}{Let $n>9$ be a nonsquare.
Then for any $n$ general points of \pr2, there are infinitely
many $m$ such that if $I(m;n)$ and $I(m;n)'$ are unhindered, then
$I(m;n)$ is rank minimal.}

\Prf We will use a criterion developed in the proof of Corollary 5.9
of [\ir{refigp}], which we briefly recall.
By \ir{igpLemma}, if $I(m;n)$ is unhindered
and $q(m;n)=0$, then $I(m;n)$ is rank minimal.
If $I(m;n)$ is unhindered, $\alpha(m;n)$ is the solution $x$ 
to the pair of inequalities
$(x+1)(x+2)-n(m+1)m>0$ and $(x+1)x-n(m+1)m\le 0$.
If in addition $I(m;n)'$ is unhindered, then
$q(m;n)=0$ if and only if also $(x+1)(x+2)-n(m+1)m-2(m+1)\le0$.
But $(x+1)(x+2)-n(m+1)m>0$ implies that $x\ge m$,
and since adding $2(x-m)$ to $(x+1)x-n(m+1)m$ gives
$(x+1)(x+2)-n(m+1)m-2(m+1)$, we see,
assuming $I(m;n)$ and $I(m;n)'$ are unhindered, that 
$q(m;n)=0$ if $(x+1)(x+2)-n(m+1)m>0$ and 
$(x+1)(x+2)-n(m+1)m-2(m+1)\le0$ have a simultaneous 
positive integer solution $x$. Let $0<\epsilon<1$; by easy but tedious 
arithmetic, one can check each real $x$ in the interval 
$[\sqrt{n}m+(\sqrt{n}-3)/2, \sqrt{n}m+(\sqrt{n}-3)/2 + \epsilon/\sqrt{n}]$
is, for $m$ sufficiently large (how large 
depending on $n$ and $\epsilon$), a solution to 
our pair of inequalities.
Thus this interval containing an integer is a criterion for
$q(m;n)$ to vanish. But by simplifying, 
there being an integer $x$ in this interval
is equivalent to there being an integer $\eta=x+1$
satisfying $0<(2\eta+1)/(2m+1)-\sqrt{n}\le 2\epsilon/((2m+1)\sqrt{n})$.

So now we consider $0<(2\eta+1)/(2m+1)-\sqrt{n}\le 
2\epsilon/((2m+1)\sqrt{n})$.
Since $n$ is not a square, we can write $n=a^2+b$
with $0<b\le 2a$, by taking $a=[\sqrt{n}]$ to be the integer part
of $\sqrt{n}$. We now show that there are infinitely many pairs of
odd integers $p,q$ such that $0<p/q-\sqrt{a^2+b}
\le 2\epsilon/(q\sqrt{a^2+b})$,
thus completing the proof.

Let $f$ and $g$ be positive odd integers
such that $f^2-(a^2+b)g^2$ is positive.
As is well known, Pell's equation, $z^2-(a^2+b)y^2=1$,
has a solution $z=c$, $y=d$ in positive integers,
and we obtain additional solutions $z=u'$, $y=v'$
from $u'+v'\sqrt{a^2+b}=(c+d\sqrt{a^2+b})^t$.
Moreover, whenever $t$ is even it is easy to check that
$u'$ is odd and $v'$ is even. Now, taking
$u+v\sqrt{a^2+b}=(u'+v'\sqrt{a^2+b})(f+g\sqrt{a^2+b})$
we obtain a solution $z=u$, $y=v$ to $z^2-(a^2+b)y^2=f^2-(a^2+b)g^2$
with $u$ and $v$ both odd. It follows that there are infinitely
many such solutions. Moreover, 
$u/v -\sqrt{a^2+b}=(f^2-(a^2+b)g^2)/(v^2(u/v+\sqrt{a^2+b}))   
<(f^2-(a^2+b)g^2)/(v^2\sqrt{a^2+b})$ which is clearly 
less than or equal to $2\epsilon/(v\sqrt{a^2+b})$ for $v$ sufficiently large. \qed

\irrnSection{When $n$ is a square}{sqrs}
Corollary 5.8 of [\ir{refigp}] shows that \ir{nnsqr} also holds for odd squares
$n>9$. In this section, we will strengthen this result for odd squares.
However, our main result in this section is \ir{evensqr}(a),
which by \ir{hlbfnctn}(i) is a precise formulation of the fact
that for each sufficiently
large $m$, if $I(m;n)$ is unhindered then it is rank minimal,
when $n$ is an even square bigger than 9. 

\prclm{Theorem}{evensqr}{Let $n=r^2$ for $r\ge 3$.
\item{(a)} Assume $r$ is even.
Then $I(m;n)$ is rank minimal for $m\ge (r-2)/4$ if 
$\alpha(m;n)=rm+r/2-1$.
\item{(b)} Assume $r$ is odd.
Then $I(m;n)$ is rank minimal for $m\ge (r^2-9)/8$ if 
$I(m;n)$ and $I(m;n)'$ are unhindered.}

The proof is at the end of this section. 

\rem{Example}{egs} The hypothesis of \ir{evensqr}(a), 
that $\alpha(m;n)=rm+r/2-1$, is in several cases known to hold.
It follows from [\ir{refNone}] and 
\ir{hlbfnctn}(i) that $\alpha(m;n)=rm+r/2-1$
for all $m>0$ when $r=4$. Thus URC holds by \ir{evensqr} 
for all $m\ge 0$ if $n=16$.
Additionally, [\ir{refCMa}, \ir{refCMb}] show that $I(m;n)$
is unhindered for any $n>9$ general points
when $m\le 12$. Consequently, by \ir{evensqr}(a) and \ir{hlbfnctn}(i), URC holds
for $I(m;r^2)$ with $12\ge m\ge (r-2)/4$ when $r>3$ 
is even. In particular: for any even square $16< n\le 50^2=2500$, 
URC holds for $n$ general points taken with multiplicity 12;
for any even square $16< n\le 46^2=2116$,
URC holds for $n$ general points taken with multiplicity 11; 
$\ldots$; and for any even square $16< n\le 14^2=196$,
URC holds for $n$ general points taken with multiplicity 3.
(For $m=2$ [\ir{refId}] shows that URC holds for 
all $n>9$.)
\vskip6pt

Our proof of \ir{evensqr} uses geometrical arguments.
Thus we now will work on the blow-up $X$ of \pr2 at the $n$ points 
$p_1,\ldots,p_n$. Results in this geometrical setting directly translate 
back to the algebraic setting to which we have up to now
mostly confined ourselves, and we refer the reader to, for example,
section 3 of [\ir{refigp}], for the dictionary to do so. 
The reader will recall from \ir{qunifunhin} and \ir{algor}
that $E_i$ is the exceptional 
locus of the blow up of $p_i$, and $L$ is the total transform
of a line in \pr2. Then the divisor classes $[L], [E_1],
\ldots, [E_n]$ give a basis of the divisor class group of $X$.
We will denote the divisor $tL-m(E_1+\cdots+E_n)$
by $F_{t,m}$, and the corresponding line bundle ${\cal O}_X(F_{t,m})$
by $\C F_{t,m}$. Note that $[L]=[F_{1,0}]$ and $[F_{t,l}]+[F_{s,m}]=[F_{t+s,l+m}]$.

\prclm{Lemma}{hlbfnctn}{Let $X$ be the blow up of \pr2 at $n$ general
points where $n=r^2$, and let $t$ and $m$ be nonnegative integers.
\item{(i)} Let $r$ be even. 
\itemitem{(a)} Then $h^1(X,\C F_{t,m})=0$ and $h^0(X,\C F_{t,m})>0$
for all $t\ge rm+(r-2)/2$.
\itemitem{(b)} If $m\ge (r-2)/4$ and $0\le t< rm+(r-2)/2$, then 
$h^0(X,\C F_{t,m})-h^1(X,\C F_{t,m})\le 0$.
\item{(ii)} Let $r$ be odd. 
\itemitem{(a)} Then $h^1(X,\C F_{t,m})=0$ and $h^0(X,\C F_{t,m})>0$
for all $t\ge rm+(r-3)/2$.
\itemitem{(b)} If $m\ge (r-1)(r-3)/(8r)$ and $0\le t<rm+(r-3)/2$, then 
$h^0(X,\C F_{t,m})-h^1(X,\C F_{t,m})\le0$.
}

\Prf Note for $t\ge0$, that $h^2(X,\C F_{t,m})=
h^0(X,\C F_{-3-t,1-m})$ by duality, and $h^0(X,\C F_{-3-t,1-m})$ 
vanishes, since $-3-t<0$. Thus for both (i)(b) and (ii)(b) we have
$h^0(X,\C F_{t,m})-h^1(X,\C F_{t,m}) = (F_{t,m}^2-K_X\cdot F_{t,m} +2)/2$
by Riemann-Roch. Now $F_{t,m}^2-K_X\cdot F_{t,m} +2=
(t+2)(t+1)-r^2m(m+1)$ is an increasing function of $t$ for $t\ge 0$,
so to prove (i)(b) and (ii)(b) 
it is enough to take $t=rm+(r-2)/2-1$ when $r$ is even and 
$t=rm+(r-3)/2-1$ when $r$ is odd. For $t=rm+(r-2)/2-1$, 
$(t+2)(t+1)-r^2m(m+1)$ becomes $(r/4)(r-4m-2)$, which
is nonpositive for $m\ge (r-2)/4$. This proves (i)(b).
For $t=rm+(r-3)/2-1$, 
$(t+2)(t+1)-r^2m(m+1)$ becomes $-2rm+(r-3)(r-1)/4$, which
is nonpositive for $m\ge (r-3)(r-1)/(8r)$. This proves (ii)(b).

To prove both parts (a), consider a specialization of
the $r^2$ points to the case of general points $p_i$ 
on a smooth plane curve $C'$ of degree $r$, and let $C$ be the proper 
transform to $X$ of $C'$. 
In the case that $r$ is even (so $t\ge rm+(r-2)/2$), the restriction
$\C F_{t,m}\otimes \C O_C$ to $C$ has degree at least $(r-2)/2+g$,
where $g$ is the genus of $C$. 
But $\C F_{t,m}\otimes \C O_C$
is a general bundle of its degree (the points $p_i$
being general points on $C'$), so
$\C F_{t,m}^{-1}\otimes \C K_C$ is a general bundle of degree
at most $g-r/2-1$, and thus has no nontrivial global sections (since 
$\hbox{Pic}^0(C)$ has dimension $g$, there are more line bundles than there
are effective divisors for any given degree less than $g$). Thus 
$h^1(C,\C F_{t,m}\otimes \C O_C)=0$ by duality. Since (i)(a) 
is true for $m=0$, our result follows for all $m\ge 0$ by induction
by taking global sections of the sequence
$$0\to \C F_{t-r,m-1}\to \C F_{t,m}\to \C F_{t,m}\otimes \C O_C\to 0$$
obtained by restriction. The result for general points of \pr2
(rather than general points of $C'$) now follows by semicontinuity.
Case (ii)(a), that $r$ is odd, is similar, except now $t\ge rm+(r-3)/2$
and $\C F_{t,m}\otimes \C O_C$ has degree at least $g-1$.
\qed

\rem{Remark}{semicont} To justify our appeal to semicontinuity
in the preceding proof, we can appeal to flat families,
using results of [\ir{refsundance}]. Alternatively, consider any 
$n$ nonnegative multiplicities $m_1,\ldots,m_n$ and 
any $n$ distinct points $p_1,\ldots,p_n\in\pr2$, with their corresponding
ideals $P_i$. We may assume that none of the ideals $P_i\subset k[x,y,z]$ contain $z$.
Let $I=\cap_{i=1}^n P_i^{m_i}$ and let $V_i$ be the vector space span
in $k[x,y]$ of the monomials of degree less than $m_i$. 
For each degree $t$ and each point $p_i$ we have a linear
map $\lambda_{ti} : R_t\to V_i$, in which a homogeneous polynomial
$f(x,y,z)$ of degree $t$ is first sent to $f(x+x_1,y+y_1,1)$,
where $(x_1,y_1)$ are affine coordinates for $p_i$ (taking $z=0$ to 
be the line at infinity), and then $f(x+x_1,y+y_1,1)$ is truncated to drop all
terms of degree $m_i$ or more. With $V$ taken to be $V_1\times\cdots\times V_n$, 
define $\Lambda_t :R_t\to V$ to be the map
$\lambda_{t1}\times\cdots\times\lambda_{tn}$.
Thus the kernel of $\Lambda_t$ is $I_t$, and the entries for the matrix
defining $\Lambda_t$ in terms of bases of monomials
are polynomials in the coordinates of the points $p_i$.
If $X$ is the blow up of \pr2 at the points $p_i$, and if we set
$\C F=\C O_X(F)$ where $F=tL-m_1E_1-\cdots-m_nE_n$, then $h^0(X,\C F)=\hbox{dim }I_t$
and by Riemann-Roch $h^0(X,\C F)-h^1(X,\C F)=\hbox{dim }R_t - \hbox{dim }V$.
Since $\hbox{dim }I_t =\hbox{dim }R_t - l$ where $l$ is the rank of
$\Lambda_t$, we see that $h^1(X,\C F) =\hbox{dim }V - l$.
But for any given $l$, 
$\Lambda_t$ having rank at least $l$ is an open condition
on the set of all $n$-tuples $(p_1,\ldots,p_n)$
of distinct points with none at infinity,
and hence $h^1(X,\C F)$ is also semicontinuous.
\vskip6pt

Our proof of \ir{evensqr}(a) involves examining certain 
specializations of the $r^2$ points.
The first specialization we will consider is easiest to specify
as a subset of the affine plane ${\bf A}^2$, rather than of \pr2.
Regard ${\bf A}^2$ as $k^2$, where $k$ is the ground field.
Consider $r$ distinct vertical lines $V_1,\ldots,V_r$ in
$k^2$ and $s=3r/2-1$ distinct horizontal lines $H_1,\ldots,H_s$.
Let $p_{ij}$, for each $1\le i\le r$ and $1\le j\le s$, be the point  
of intersection of $V_i$ with $H_j$. We choose our $r^2$ points 
in blocks from among these
$rs$ points $p_{ij}$. The first block is 
$B_1=\{p_{ij} : 1\le i\le2, 1\le j\le s\}$, the second is
$B_2=\{p_{ij} : 3\le i\le4, 1\le j\le s-2\}$, etc., and the last block is
$B_{r/2}=\{p_{ij} : r-1\le i\le r, 1\le j\le s-(r-2)\}$.
Note that all together, the union of these $r/2$ blocks
contains $2s+\cdots+2(s-(r-2))=2(sr/2-(2+4+\cdots+(r-2))
=sr-4(1+\cdots+(r/2-1))=(3r/2-1)r-4(r/2-1)(r/2)/2=r^2$
points. It will be convenient to denote $rm+(r-2)/2$ by $t_m$
and by $\mu_{t_m}$ the map $H^0(X, \C F_{t_m,m})\otimes H^0(X,\C F_{1,0})\to
H^0(X, \C F_{t_m+1,m})$ given by multiplication. 

\prclm{Proposition}{blcks}{With respect to the configuration
$B_1\cup\cdots\cup B_{r/2}$ of $r^2$ points of \pr2 specified
above, $\mu_{t_1}$ is surjective and $h^1(X, \C F_{t_1,1})=0$.}

\Prf Apply the method and results of [\refGGR]. Let $I\subset R$ be the homogeneous
ideal of the $r^2$ points in the homogeneous coordinate ring $R$ of \pr2. 
Let $f\in R$ be a homogeneous element of degree 1
not vanishing on any of the $r^2$ points; thus the image of $f$ 
in the quotient $R/I$ is not a zero divisor.
Hence, as discussed in [\refGGR], $J=I+(f)/(f)$ has the same minimal
number of homogeneous generators
in every degree as does $I$, but $J$ is a monomial ideal which, by the discussion
in [\refGGR], is easy to handle explicitly. The result in our case is that
there are no generators in degrees greater than $3r/2-1$, hence the same is true
for $I$, which shows that $\mu_{t_1}$ is surjective, as required. 
To see that $h^1(X, \C F_{t_1,1})=0$, it is enough to check that
the points impose independent conditions on forms of degree $t_1$;
i.e., that $h^0(X, \C F_{t_1,1})=(t_1+2)(t_1+1)/2-r^2$.
But the fact that $J=I+(f)/(f)$ has the same number of generators in every 
degree as does $I$ means that $I$ has one generator in degree $r$
and $r/2$ in degree $t_1=3r/2-1$. Thus the dimension of $I$ in
degree $t_1$ is $(r/2) + ((r/2)(r/2+1)/2)$, which is indeed
$(t_1+2)(t_1+1)/2-r^2$.\qed

\prclm{Proposition}{rdctntcrv}{Let 
$r\ge0$ be even with $m\ge0$. Then for general points $p_1,\ldots,p_{r^2}$
of a general smooth plane curve $C'$ of degree $r$, 
taking $X$ to be the blow up of
\pr2 at $p_1,\ldots,p_{r^2}$, the map
$\mu_{t_m}$ is surjective and $h^1(X, \C F_{t_m,m})=0$
for all $m\ge 0$.}

\Prf First let $C'$ be any smooth plane curve of degree $r$ 
and let $C$ be the proper transform of $C'$ to $X$.
As in the proof of \ir{hlbfnctn}, we have the exact sequence 
$0\to \C F_{t_m-r,m-1}\to \C F_{t_m,m}\to \C F_{t_m,m}\otimes \C O_C\to 0$,
and $h^1(X, \C F_{t_m,m})=0$ for all $m\ge0$.
The exact sequence leads by the snake lemma (see [\refMu] or
[\reffatpts]) to an exact sequence\par
\vskip6pt
\noindent\hbox to\hsize{\hbox to0in{$(*)$\hss}\hfil
$\hbox{cok }(\mu_{t_m-r}) \to \hbox{cok }(\mu_{t_m})\to 
\hbox{cok }(\mu_{C,t_m})\to 0,$\hfil}
\vskip6pt
\noindent
where $\mu_{C,t_m}$ is the map $H^0(C, \C F_{t_m,m}\otimes 
\C O_C)\otimes H^0(X,\C F_{1,0})\to
H^0(X, \C F_{t_m+1,m}\otimes \C O_C)$. But
$\C F_{t_m,m}\otimes \C O_C$ is a general bundle of degree 
$r(r-2)/2$, regardless of $m$, so if $\hbox{cok }(\mu_{C,t_1})=0$,
then $\hbox{cok }(\mu_{C,t_m})=0$ for all $m\ge 1$, and
$\hbox{cok }(\mu_{C,t_1})=0$ follows from $(*)$
if we show that $\hbox{cok }(\mu_{t_1})=0$. But $\hbox{cok }(\mu_{t_1})=0$
and $h^1(X, \C F_{t_1,1})=0$
for the configuration of points given in \ir{blcks}, and these
points are points of a plane curve $C'$ of degree $r$
(the union of $r$ lines),
so by semicontinuity (see \ir{semicontB}) $\hbox{cok }(\mu_{t_1})=0$
holds for general points of a general curve $C'$ of degree $r$.

Finally, induction using $(*)$ gives $\hbox{cok }(\mu_{t_m})=0$ for all 
$m\ge 0$. (Note that $\mu_{t_m-r}=\mu_{t_{m-1}}$ so 
$\mu_{t_m-rm}=\mu_{t_{0}}$, and 
$\hbox{cok }(\mu_{t_m-rm})=0$ since
$\C F_{t_m-rm,0}=\C F_{t_0,0}$ can be regarded 
as $\C O_{\pr2}(((r-2)/2)L)$ on \pr2, where the result is obvious.)
\qed

\rem{Remark}{semicontB} We now show that the requirements
that $\hbox{cok }(\mu_{t_1})=0$ and $h^1(X, \C F_{t_1,1})=0$
together impose an open condition on 
$r^2$-tuples of points $(p_1,\ldots,p_{r^2})$.
More generally, let $I$ and \C F be as in \ir{semicont}. As there, we have 
the map $\Lambda_t: R_t\to V$.
Then $I_t\otimes R_1$ is the kernel of 
$\Lambda_t\otimes \hbox{id}_{R_1}: R_t\otimes R_1\to V\otimes R_1$
and the kernel of $\mu_t: I_t\otimes R_1\to I_{t+1}$
is also the kernel of
$\gamma: R_t\otimes R_1\to (V\otimes R_1)\oplus R_{t+1}$,
where $\gamma$ is $(\Lambda_t\otimes \hbox{id}_{R_1})\oplus \mu_t'$,
and where $\mu_t: I_t\otimes R_1\to I_{t+1}$ and 
$\mu_t': R_t\otimes R_1\to R_{t+1}$ are given by multiplication.
(This is consistent with our usage above, since sections of line bundles
on blow ups of \pr2 can be identified with subspaces of $R_t$ for 
appropriate $t$.) Since, as in \ir{semicont}, $\gamma$ can be given
by a matrix whose entries are polynomials in the coordinates
of the points $p_i$, the rank of $\gamma$ and hence the 
dimension of $\hbox{ker }(\mu_t)$ is semicontinuous.
But $\hbox{dim cok }(\mu_t)=h^0(X, \C F(L))-3h^0(X, \C F)
+\hbox{dim ker }(\mu_t)$, which is 
$\hbox{dim cok }(\mu_t)=t+F\cdot K_X-F^2
+\hbox{dim ker }(\mu_t)$ if $h^1(X, \C F)=0$,
and since $h^1(X, \C F)$ is also semicontinuous, it
follows that it is an open condition to
require that both $\hbox{cok }(\mu_t)$ 
and $h^1(X, \C F)$ vanish.
\vskip6pt

\vskip6pt

\noindent{\bf Proof} (of \ir{evensqr}): (a) Assume 
that $\alpha(m;n)=rm+r/2-1$.
Then by definition of $\alpha(m;n)$, we have $I(m;n)_t=0$ for $t<rm+r/2-1$
and by \ir{hlbfnctn}(i)(a) we have $h^1(X, \C F_{t,m})=0$
for $t\ge rm+r/2-1$. Thus $I(m;n)$ is unhindered.
Moreover, since $I(m;n)_t=0$ for $t<rm+r/2-1$, there are no homogeneous
generators for $I(m;n)$ in degrees less than $rm+r/2-1$
and there are $h_{I(m;n)}(\alpha(m;n))$ generators
in degree $rm+r/2-1$. Clearly, then, $I(m;n)$ is rank minimal
if $F_0$ has rank $h_{I(m;n)}(\alpha(m;n))$, so it is enough 
to prove that these $h_{I(m;n)}(\alpha(m;n))$ elements
of degree $rm+r/2-1$ generate $I(m;n)$.

By \ir{hlbfnctn}(i), the
Castelnuovo-Mumford regularity of
$I(m;n)$ is at most $\alpha(m;n)+1$, so no generators 
need be taken in degrees greater than $rm+r/2$
(see [\ir{refDGM}] or Lemma 2.9 of [\reffatpts]).
Hence we now need only show that no generators
need be taken in degree $rm+r/2$; i.e., 
that $\mu_{t_m}: H^0(X, \C F_{t_m,m})\otimes H^0(X,\C F_{1,0})\to
H^0(X, \C F_{t_m+1,m})$ is surjective,
where, as in \ir{rdctntcrv}, $t_m=rm+r/2-1$.
But this follows by \ir{rdctntcrv} and 
semicontinuity (cf. \ir{semicontB}).

(b) (We note that when $n=r^2>9$ is an odd square,
[\ir{refigp}] noted but did not explicitly show that 
$I(m;r^2)$ is rank minimal for all but finitely many $m$
for which $I(m;r^2)$ and $I(m;r^2)'$ are unhindered.)
By \ir{igpLemma}(b), if $I(m;r^2)$ is unhindered
and $q(m;r^2)=0$, then $I(m;r^2)$ is rank minimal.
But \ir{hlbfnctn}(ii) implies that
$\alpha(m;r^2)=rm+(r-3)/2$ for $m\ge (r-1)(r-3)/(8r)$.
Now, $q(m;r^2)=h_{I(m;n)'}(\alpha(m;r^2))$, and
using $t=\alpha(m;r^2)$ and assuming that
$I(m;r^2)'$ is unhindered, we have
$h_{I(m;n)'}(\alpha(m;r^2))=
\hbox{max}\{0,(t+2)(t+1)/2-(r^2-1)m(m+1)/2-(m+1)(m+2)/2\}$.
But $(t+2)(t+1)/2-(r^2-1)m(m+1)/2-(m+1)(m+2)/2\le 0$
for $r\ge3$ with $m\ge (r^2-9)/8$, and since
$(r^2-9)/8\ge (r-1)(r-3)/(8r)$ when $r\ge 3$,
we see that $I(m;r^2)$ is rank minimal
if $m\ge (r^2-9)/8$ for $r\ge3$, whenever 
$I(m;r^2)$ and $I(m;r^2)'$ are unhindered.
\qed

\irrnSection{Additional Results}{addres}
Our results on resolutions of uniform ideals in certain cases 
extend to quasiuniform ideals. This raises the question
of whether quasiuniform ideals involving general points
of \pr2 are always rank minimal; we know of no counterexamples.

\prclm{Proposition}{qunifresprop}{Assume $q(m;r)=0$ 
for some $m$ and some $r\ge 2$. 
Let $n\ge r$ and let ${\bf m}=(m_1,\ldots,m_n)$ where $m=m_1=\cdots=m_r$
and $\sum_{i>r} (m_i^2+m_i)/2<h_{I(m;r)}(\alpha(m;r))$.
If $I({\bf m};n)$ is unhindered, then 
$I({\bf m};n)$ is rank minimal.}

\Prf By the hypothesis $\sum_{i>r} (m_i^2+m_i)/2<h_{I(m;r)}(\alpha(m;r))$,
we see $\alpha(m;r)=\alpha({\bf m};n)$, 
so $q(m;r)=0$ implies $q({\bf m};n)=0$, and the result follows
by \ir{igpLemma}(b). \qed

Note that if $r\ge 9$ in \ir{qunifresprop}, 
then, as conjectured by QHC, we expect that  
$I({\bf m};n)$ is indeed unhindered if ${\bf m}$ is quasiuniform.
Similarly, we expect that the hypotheses of unhinderedness
in the following result always hold (since
after reordering $m_i$, if need be, 
for $i>r^2$, ${\bf m}$ is quasiuniform).

\prclm{Corollary}{lastcor}{Let $n\ge r^2$ where $r\ge3$ is odd 
and let $m=m_1=\cdots=m_{r^2}\ge (r^2-9)/8\ge
\sum_{r^2<i\le n} (m_i^2+m_i)/2$.
If, for ${\bf m}=(m_1,\ldots,m_n)$,  
$I({\bf m};n)$, $I(m;r^2)$ and $I(m;r^2)'$ are unhindered,
then $I({\bf m};n)$ is rank minimal.}

\Prf As in the proof of \ir{evensqr}(b),
$I(m;r^2)$ and $I(m;r^2)'$ being unhindered
implies $q(m;r^2)=0$. From \ir{hlbfnctn}(ii) we find 
$h_{I(m;r^2)}(\alpha(m;r^2))=(r^2-1)/8$,
so \ir{qunifresprop} gives the result.
\qed

We close with a result involving no conditional hypotheses
of unhinderedness.

\prclm{Proposition}{lastprop}{Let $m>0$, let $t\ge 0$,
and choose $n\ge9$ such that
$-m-1\le n - [9+3tm+(t+1)(t+2)/2]\le m$.
Then $I({\bf m};n)$ is rank minimal for ${\bf m}=(m_1,\ldots,m_n)$ with
$m=m_1=\cdots=m_{9}$ and $m_i=1$ for $9<i\le n$.}

\Prf It is known (see [\ir{refanti}]) that 
\ir{genconj} holds for 9 general points,
and thus (by \ir{qunifunhin}) so does QHC, but general simple points
impose independent conditions, so
QHC also holds when ${\bf m}$ satisfies the hypotheses 
of \ir{lastprop}. Thus
$I({\bf m};n)$, $I({\bf m};n)'$ and $I({\bf m};n)''$
are unhindered.

Given this, if $8+(3t-1)m+(t+1)(t+2)/2\le n\le 
8+3tm+(t+1)(t+2)/2$, we can easily check that
$\alpha({\bf m};n)=3m+t$ and then that 
$q({\bf m};n)=0$, whereas if
$9+3tm+(t+1)(t+2)/2\le n\le 
9+(3t+1)m+(t+1)(t+2)/2$, then
$\alpha({\bf m};n)=3m+t+1$ and 
$l({\bf m};n)>0$. Now the result
follows by \ir{igpLemma}.
\qed

\References

\bibitem{\refCam} Campanella, G. {\it Standard 
bases of perfect homogeneous
polynomial ideals of height $2$}, 
J.\ Alg.\  101 (1986), 47--60.

\bibitem{\refCat} Catalisano, M.\ V. {\it ``Fat'' points on a conic}, 
Comm.\ Alg.\  19(8) (1991), 2153--2168.

\bibitem{\refCMa} Ciliberto, C. and Miranda, R. {\it Degenerations
of planar linear systems}, 
J. Reine Angew. Math. 501 (1998), 191-220.

\bibitem{\refCMb} Ciliberto, C. and Miranda, R. {\it Linear systems
of plane curves with base points of equal multiplicity}, 
Trans AMS, to appear.

\bibitem{\refDGM} Davis, E.\ D., Geramita, A.\ V., and 
Maroscia, P. {\it Perfect
Homogeneous Ideals: Dubreil's Theorems Revisited},
Bull.\ Sc.\ math., $2^e$ s\'erie, 108 (1984), 143--185.

\bibitem{\refDu} Dubreil, P.
{\it Sur quelques propri\'et\'es des syst\`emes de points dans
le plan et des courbes gauches alg\'ebriques},
Bull.\ Soc.\ Math.\ France, 61 (1933), 258--283.

\bibitem{\refFione} Fitchett, S. {\it Doctoral dissertation}, 
University of Nebraska-Lincoln, 1996.

\bibitem{\refFitwo} Fitchett, S. {\it Maps of linear systems
on blow ups of the projective plane}, 
preprint, 1997.

\bibitem{\refFHH} Fitchett, S., Harbourne, B., and Holay, S.
{\it Resolutions of Fat Point Ideals for $n\le8$ General Points of \pr2}, 
in preparation, 1999.

\bibitem{\refGGR} Geramita, A.\ V., Gregory, D.\ and Roberts, L.
{\it Minimal ideals and points in projective space},
J.\ Pure and Appl.\ Alg. 40 (1986), 33--62.

\bibitem{\refGi} Gimigliano, A. {\it Our thin knowledge of fat points},
Queen's papers in Pure and Applied Mathematics, no. 83,
The Curves Seminar at Queen's, vol. VI (1989).

\bibitem{\refvanc} Harbourne, B. {\it The geometry of 
rational surfaces and Hilbert
functions of points in the plane},
Can.\ Math.\ Soc.\ Conf.\ Proc.\ 6 
(1986), 95--111.

\bibitem{\refsundance}  \manyby. {\it Iterated blow-ups and moduli for 
rational surfaces}, in: Algebraic Geometry, 
Sundance 1986, LNM \#1311 (1988), 101--117.

\bibitem{\refravello} \manyby. {\it 
Points in Good Position in \pr 2}, in:
Zero-dimensional schemes, Proceedings of the
International Conference held in Ravello, Italy, June 8--13, 1992,
De Gruyter, 1994.

\bibitem{\refanti} \manyby.  {\it Anticanonical 
rational surfaces}, Trans. A.M.S. 349 (1997), 1191--1208.

\bibitem{\reffatpts} \manyby. {\it Free Resolutions of Fat Point 
Ideals on \pr2}, JPAA 125 (1998), 213--234.

\bibitem{\refigp} \manyby. {\it The Ideal Generation 
Problem for Fat Points}, preprint, to appear, JPAA.

\bibitem{\refHi} Hirschowitz, A.
{\it Une conjecture pour la cohomologie 
des diviseurs sur les surfaces rationelles g\'en\'eriques},
Journ.\ Reine Angew.\ Math. 397
(1989), 208--213.

\bibitem{\refId} Ida, M.
{\it The minimal free resolution for the first infinitesimal 
neighborhoods of $n$ general points in the plane},
J.\ Alg., to appear.

\bibitem{\refMir} Miranda, R. {\it Linear Systems of Plane
Curves}, Notices AMS 46(2) (1999), 192--202.

\bibitem{\refMu} Mumford, D. {\it Varieties defined by quadratic equations},
in: Questions on algebraic varieties, Corso C.I.M.E. 1969 Rome: Cremonese,
1969, 30--100.

\bibitem{\refNone} Nagata, M. {\it On the 14-th problem of Hilbert}, 
Amer.\ J.\ Math.\ 33 (1959), 766--772.

\bibitem{\refNtwo} \manyby. {\it On rational surfaces, II}, 
Mem.\ Coll.\ Sci.\ 
Univ.\ Kyoto, Ser.\ A Math.\ 33 (1960), 271--293.

\bibitem{\refR} Ro\'e, J. {\it On the existence of plane curves
with imposed multiple points}, 
preprint (1998).

\bye